\documentclass{article}

\usepackage{graphicx} 

\usepackage{latexsym, amscd, amsfonts, eucal, mathrsfs, amsmath, amssymb, amsthm, xypic,xr, makeidx, stmaryrd}

\usepackage{tikz-cd}
\usepackage{bookmark}
\usepackage{makeidx}
\usepackage{geometry}
\usepackage{indentfirst}
\usepackage{enumitem}
\usepackage{mathrsfs}
\usepackage{enumitem}
\usepackage{verbatim}
\usepackage{upgreek}
\usepackage{textcomp}
\usepackage{amsmath}
\usepackage{dsfont}

\newtheorem{thm}{Theorem}[subsection]
\newtheorem{defn}[thm]{Definition}
\newtheorem{notation}[thm]{Notation}

\newtheorem{example}[thm]{Example}

\AfterEndEnvironment{thm}{\noindent\ignorespaces}
\AfterEndEnvironment{defn}{\noindent\ignorespaces}
\AfterEndEnvironment{notation}{\noindent\ignorespaces}
\AfterEndEnvironment{prop}{\noindent\ignorespaces}
\AfterEndEnvironment{lem}{\noindent\ignorespaces}
\AfterEndEnvironment{cor}{\noindent\ignorespaces}
\AfterEndEnvironment{example}{\noindent\ignorespaces}

\DeclareMathOperator{\sym}{\mathbb{S}}
\DeclareMathOperator{\X}{{\mathcal{X}}}
\DeclareMathOperator{\OX}{\mathcal{O}_{\X}}
\DeclareMathOperator{\G}{{\mathcal{G}_{\sym} ^{alg}}}
\DeclareMathOperator{\T}{{\mathcal{T}_{\sym} ^{alg}}}
\DeclareMathOperator{\shv}{\mathrm{Shv}}
\DeclareMathOperator{\spec}{\mathrm{Spec}}
\DeclareMathOperator{\fun}{\mathrm{Fun}}
\DeclareMathOperator{\map}{\mathrm{Map}}
\DeclareMathOperator{\pro}{\mathrm{Pro}}
\DeclareMathOperator{\ind}{\mathrm{Ind}}
\DeclareMathOperator{\perf}{\mathrm{Perf}}

\title{Some memos on Stable Symplectic Structured Space}
\author{Eita Haibara, Shelly Miyano}
\date{April 2025}

\begin{document}
\maketitle

\begin{abstract}
In these memos, we define a pregeometry $\T$ and a geometry $\G$ which integrate symplectic manifolds with $E_{\infty}$-ring sheaves, enabling the construction of $\G$-schemes as structured $\infty$-topoi. Our framework and results establish a profound connection between algebraic invariants and homological properties, opening new pathways for exploring symplectic phenomena through the lens of higher category theory and derived geometry.
\end{abstract}

\tableofcontents

\section{Background}
This is just a background section, so we'll skip the detailed proof. We assume that the reader reads \cite{LuDAGV}, \cite{kitchloo12} ,\cite{ha}, and \cite{costello}.
\newline
\subsection{Calabi-Yau}
A d-Calabi-Yau $\infty$-category is a stable $\infty$-category over a field k equipped with a Calabi-Yau structure, meaning its Serre functor $S_C$ is isomorphic to the shift functor [d], where d is an integer called the dimension. For a dg-category or stable $\infty$-category, this condition implies a duality: for objects X,Y $\in$ C, there is a natural isomorphism
\begin{center}
$Hom_C(X,Y) \cong Hom_C(Y,X[d])^*$
\end{center}
where * denotes the dual over the base field k. A Calabi-Yau structure typically includes a non-degenerate pairing compatible with this duality.

\section{Pregeometry $\T$ and Geometry $\G$} 
In this section, I'm sharing the results of my research on Symplectic geometry that I did during my free time for a week. But I don't know if I'm going to further this research. I hope someone will continue my research.
\begin{flushright}
    by Shelly Miyano.
\end{flushright}
\begin{notation}
We let $M^{-\mathcal{\uptau}}$ denote the Thom spectrum associated with the stable normal bundle of the symplectic manifold and $\mathbb{S}$ denote the stable symplectic category.\cite{kitchloo12}
\end{notation}

 Note that:$M^{-\mathcal{\uptau}}$ is the $E_\infty$-ring spectrum according to Remark 7.5 of \cite{kitchloo12}.

\subsection{Pregeometry $\T$}
\begin{defn}
Pregeometry $\T$ consists of the following data: 
\\
\newline
(1) Objects: Pair $(M,M^{-\mathcal{\uptau}})$, where $(M,\mathcal{\upomega}_M)$ is a symplectic manifold.
\newline
\newline
(2) Morphisms: $(f, \pi): (M,M^{-\mathcal{\uptau}}) \to (N,N^{-\mathcal{\uptau}})$. $f: M \to N$ is a symplectic map satisfying $f^{*}\mathcal{\upomega}_N=\mathcal{\upomega}_M$, where $f^{*}$ denotes the pullback of the symplectic form $\mathcal{\upomega}_N$ on $N$ to $M$. $\pi : N^{-\mathcal{\uptau}} \to M^{-\mathcal{\uptau}}$ is a morphism of $E_\infty$-rings. Since $N^{-\mathcal{\uptau}}$ and $M^{-\mathcal{\uptau}}$ are Thom spectra with $E_\infty$-ring structure, $\pi$ is a map that respects this multiplicative structure, opearating contravariantly to the direction of $f$.                
\newline
\newline
(3) Admissible morphisms: $(f,\pi): : (M,M^{-\tau}) \to (N,N^{-\tau})$ os admissible if either: $f:M \to N$ is an open symplectic embedding, and $\pi: N^{-\tau} \to M^{-\tau}$ is the restriction map of Thom spectra or ($f,\pi)$ is an isomorphism. (i.e. $f$ is a symplectomorphism)
\newline
\newline
(4) Grothendieck topology: The topology on Pregeometry $\T$ is generated by collection of admissible morphisms $[(i_\alpha, \pi_\alpha): (U_\alpha,U^{-\mathcal{\uptau}}_\alpha) \to (M,M^{-\mathcal{\uptau}})]_{\alpha\in I}$, where $[i_\alpha:U_\alpha \to M]$ forms an open cover of M in the symplectic topology. 

\end{defn}
\\
Note that: (4) allows for the gluing of local data, similar to how affine schemes are patched together in algebraic geometry, but here adapted to the symplectic setting with Thom spectra. 
\\
\\
Note that: $\T$ has finite limits.
\newline
\begin{defn}
$\T$-structure on $\infty$-topos $\X$ is a functor $\mathcal{O}: \T \to \mathcal {X}$ which is left exact and satisfies the sheaf condition with respect to the Grothendieck topology on $\T$. Specifically, for any object $(M,M^{-\mathcal{\uptau}}) \in \T$ and any admissible covering family $[(U_i,U^{-\mathcal{\uptau}}_i) \to (M,M^{-\mathcal{\uptau}})]$, the natural map: 
\begin{equation*}
\mathcal{O}(M,M^{-\mathcal{\uptau}}) \to \displaystyle\lim_{\Delta} \prod \mathcal{O}(U_{i_0} \times _M \cdots \times _M U_{i_n},(U_{i_0} \times _M \cdots \times _M U_{i_n})^{-\mathcal{\uptau}})
\end{equation*}
must be an equivalence in $\X$.
\end{defn}

This ensures that $\mathcal{O}$ glues local data consistently over coverings, much like a sheaf in classical geometry.  
\\
\newline

\subsection{Geometry $\G$}

\begin{defn}
Geometry $\G$ consists of the following data:
\\
\newline
(1) Objects: Pairs $(M,A)$, where $M$ is a symplectic manifold with symplectic form $\mathcal{\upomega}_M$. $A$ is a sheaf of $E_\infty$-rings, and these assignments satisfy the sheaf condition: for an open cover $(U_i \to U)_{i \in I }$ , the diagram 

\begin{equation*}
A(U) \to \prod_{i}A(U_i) \rightrightarrows \prod_{i,j} A(U_i \cap U_j)
\end{equation*}
is an equalizer in the category of $E_\infty$-rings. This ensures $A$ glues consistently over open covers.
\newline
\newline
(2) Morphisms: A morphism $(f,\pi): (M,A) \to (N,B)$, where $f: M \to N$ is a symplectic map with $f^{*}\mathcal{\upomega}_N=\mathcal{\upomega}_M$, preserving the symplectic structure. $\pi: f^{*}B \to A$ is a morphism of sheaves of $E_\infty$-rings, where $f^{*}B$ is the pullback sheaf defined by $(f^{*}B)(U) = B(f(U))$ for open $U \subseteq M$.
\newline
\newline
(3) Admissible morphisms: A morphism ($f,\pi) : (M,A) \to (N,B)$ os admissible if either: $f:M \to N$ is an open symplectic embedding, and $\pi:f^*B \to A$ is an isomorphism of sheaves of $E_{\infty}$-rings or ($f,\pi)$ is an isomorphism. (i.e. $f$ is a symplectomorphism)
\newline
\newline
(4) Grothendieck topology: The topology on $\G$ is generated by collections $[(U_i,A|_U) \to (M,A)]_{i \in I}$, where $(U_i \to M)_{i \in I}$ is an open cover of M in the symplectic topology. 
\end{defn}

Note that (4) allows objects to be reconstructed from local data via descent, enabling the definition of schemes or stacks in this context.
\newline
\newline
Note that: The terminal object in $\G$ is ({$\star$},k), where {$\star$} is a single point, viewed as a 0-dimensional symplectic manifold with the trivial symplectic form $\omega=0$, and $k$ is the base field, considered as a constant sehaf of $E_{\infty}$-rings on {$\star$}. 

\begin{defn}
\begin{itemize}
    The product of two objects \((M, A)\) and \((N, B)\) in \(\mathcal{G}_S^{\text{alg}}\) is \((M \times N, A \boxtimes B)\), where:
        \item \(M \times N\) is the Cartesian product of the symplectic manifolds \(M\) and \(N\), equipped with the symplectic form \(\omega_{M \times N} = \pi_M^* \omega_M + \pi_N^* \omega_N\). Here, \(\pi_M: M \times N \to M\) and \(\pi_N: M \times N \to N\) are the projection maps, and \(\omega_M\), \(\omega_N\) are the symplectic forms on \(M\) and \(N\), respectively.
        \item \(A \boxtimes B = \pi_M^* A \otimes \pi_N^* B\) is the external tensor product of the sheaves \(A\) and \(B\), forming a sheaf of \(E_\infty\)-rings on \(M \times N\). The tensor product is taken over the base field \(k\).
    Projection Morphisms:
        \item \((p_M, \pi_A): (M \times N, A \boxtimes B) \to (M, A)\), where \(p_M: M \times N \to M\) is the projection, and \(\pi_A: p_M^* A \to A \boxtimes B\) is the natural inclusion into the tensor product.
        \item \((p_N, \pi_B): (M \times N, A \boxtimes B) \to (N, B)\), where \(p_N: M \times N \to N\) is the projection, and \(\pi_B: p_N^* B \to A \boxtimes B\) is the natural inclusion into the tensor product.
    Universal Property: For any object \((P, C)\) with morphisms \((f, \pi): (P, C) \to (M, A)\) and \((g, \rho): (P, C) \to (N, B)\), there exists a unique morphism \((h, \sigma): (P, C) \to (M \times N, A \boxtimes B)\) such that the following diagrams commute:
    \[
    (p_M, \pi_A) \circ (h, \sigma) = (f, \pi) \quad \text{and} \quad (p_N, \pi_B) \circ (h, \sigma) = (g, \rho).
    \]
\end{itemize}
\end{defn}

\begin{defn}
\begin{itemize}
    Given two morphisms \((f, \pi): (M, A) \to (P, C)\) and \((g, \rho): (N, B) \to (P, C)\) in \(\mathcal{G}_S^{\text{alg}}\), their pullback is \((M \times_P N, D)\), where:
   
        \item \(M \times_P N\) is the fiber product of the symplectic manifolds \(M\) and \(N\) over \(P\), i.e., 
        \[
        M \times_P N = \{(m, n) \in M \times N \mid f(m) = g(n)\},
        \]
        equipped with a symplectic form inherited from \(M\) and \(N\). This assumes that \(f\) and \(g\) satisfy conditions (e.g., transversality) ensuring \(M \times_P N\) is a symplectic manifold.
        \item \(D = A \times_C B\) is the pullback of the sheaves \(A\) and \(B\) over \(C\) in the category of sheaves of \(E_\infty\)-rings on \(M \times_P N\). Explicitly, \(D\) is the sheaf on \(M \times_P N\) whose sections over an open set \(U \subset M \times_P N\) are pairs \((a, b) \in A(U_M) \times B(U_N)\) (where \(U_M\) and \(U_N\) are the projections of \(U\) to \(M\) and \(N\), respectively) such that \(\pi(a) = \rho(b)\) in \(C(f(U_M)) = C(g(U_N))\).
    
    Projection Morphisms:
   
        \item \((p_M, \alpha): (M \times_P N, D) \to (M, A)\), where \(p_M: M \times_P N \to M\) is the projection, and \(\alpha: p_M^* A \to D\) is the map sending \(A\) to the first component of \(D\).
        \item \((p_N, \beta): (M \times_P N, D) \to (N, B)\), where \(p_N: M \times_P N \to N\) is the projection, and \(\beta: p_N^* B \to D\) is the map sending \(B\) to the second component of \(D\).
    
    Universal Property: For any object \((Q, E)\) with morphisms \((h, \sigma): (Q, E) \to (M, A)\) and \((k, \tau): (Q, E) \to (N, B)\) such that 
    \[
    (f, \pi) \circ (h, \sigma) = (g, \rho) \circ (k, \tau),
    \]
    there exists a unique morphism \((u, \phi): (Q, E) \to (M \times_P N, D)\) making the appropriate diagrams commute.
\end{itemize}
\end{defn}

\begin{defn}
$\G$-structure on $\infty$-topos $\X$ is a functor $\mathcal{O}: \G \to \X$ that is left exact and satisfies the sheaf condition with respect to the topology on $\G$: For a cover $ \left[(U_i,A|_{U_1}) \to (M,A)\right]$, the map
\begin{equation*}
\mathcal{O}(M,A) \to \varprojlim(\prod_i \mathcal{O}(U_i,A|_{U_i}) \rightrightarrows \prod_{i,j}\mathcal{O}(U_i \cap U_j, A|_{U_i \cap U_j}) \cdots)
\end{equation*}
is an equivalence in $\X$, ensuring consistent gluing over covers.

\end{defn}
\begin{defn}
A pair $(\X,\mathcal{O})$, where $\X$ is an $\infty$-topos and $\mathcal{O}$ is a $\G$-structured sheaf, is called $\G$-structured $\infty$-topos.
\end{defn}
\newline
Category of Structured $\infty$-Topoi:
\begin{notation}
We let $\mathcal{L}\mathcal{T}op(\G)$ denote the $\infty$-category of pairs ($\X$,$\OX$).
\end{notation}

\begin{example}
Consider $\X$ = S, the $\infty$-category of spaces (which is an $\infty$-topos).
\newline
\newline
Definition:
For every object $(M,A) \in \G$, define $\mathcal{O}(M,A)=*$, the one-point space. For every morphism $(f,\pi): (M,A) \to (N,B)$, define $\mathcal{O}(f,\pi):*\to *$ as the identity map.
\newline
\newline
Verification: The constant functor to a terminal object preserves all limits, since $*$ is terminal in $S$.(Left exactness) For a cover $[(U_i,A|_{U_i}) \to (M,A)]$, we have $\mathcal{O}(M,A)=*$ and $\lim \mathcal{O} (U_i,A|_{U_i}) = \lim * = *$, so the map is equivalence.(Locality)
\newline
\newline
This is a valid $\G$-structured sheaf.
\end{example}

Note that: There is a difference in the definition of (pre)geometry \cite{LuDAGV}  because of idempotent complete.

\subsection{Scheme theory}
We assume that the reader reads \cite{LuDAGV}.
\begin{example}
Let's construct a $\G$-structured sheaf associated to specific object, mimicking the construction of affine scheme. Take $(M,A) \in \G$, where $M = \mathbb{R}^2$ is equipped with standard symplectic form $\mathcal{\upomega}_M = dx\wedge dy$, and $A = C_M^\infty$ is the sheaf of smooth real valued functions on $M$ (viewed as a sheaf of commuative rings, hence $E_\infty$-rings in a discrete sense). Define $S = Hom_{\G}( -,(M,C_M^\infty)) \in \ind((\G)^{op})$ by the Yoneda embedding. Then the affine $\G$-scheme is $\spec^{\G}(S) = (\X,\OX)$, where $\X = Shv(M)$, the $\infty$-topos of sheaves on M (with the open cover topology). Since $S$ is representable by $(M,C_M^\infty)$, and the topology on $\G$ is based on open covers of $M$, we identifiy $\shv(S) \cong \shv(M)$. $\OX: \G \to \shv(M)$ is defined as follows:
\newline
For each $(N,B) \in \G$, $\OX (N,B)$ is the sheaf on $M$ that assigns to each open $U \subseteq M$:
\begin{equation*}
\OX (N,B)(U) = \map_{\G}((U,C_U^\infty),(N,B)).
\end{equation*}
For a morphism $(f, \pi): (N,B) \to (N',B')$, $\OX(f, \pi): \OX(N',B') \to \OX(N,B)$ is the natural transformation induced by post-composition with $(f, \pi)$.
\newline
We will give an example of $\mathcal{O}(N,B)$. Let's evaluate $\OX(N,B)$ for specific choice of $(N,B)$:
\newline
\newline
\underline{Case 1}: $(N,B) = (\mathbb{R}^2,C_{\mathbb{R}^2}^\infty)$
\newline
Computation: $\mathcal{O}(\mathbb{R}^2,C_{\mathbb{R}^2}^\infty)(U) = \map_{\G}((U,C_U^\infty),(\mathbb{R}^2,C_{\mathbb{R}^2}^\infty))$
\newline
Morphisms: A pair $(f, \pi)$, where: $f: U \to \mathbb{R}^2$ is a symplectic map. $\pi: f^{*}C_{\mathbb{R}^2}^\infty \to C_U^\infty$ assigns to each smooth function on f(U) a smooth function on U. Since f is a diffeomorphism, $f^{*}C_{\mathbb{R}^2}^\infty \cong C_U^\infty$, and $\pi$ can be the natural isomorphism.
\newline
\newline
\underline{Case 2}: (N,B) = $(\mathbb{R}^4,C_{\mathbb{R}^4}^\infty)$
\newline
Computation: $\OX(\mathbb{R}^4,C_{\mathbb{R}^4}^\infty)$(U) = $\map_{\G}((U,C_U^\infty),(\mathbb{R}^4,C_{\mathbb{R}^4}^\infty))$
\newline
Morphisms: A pair $(f, \pi)$, where: $f: U \to \mathbb{R}^4$ satisfies $f^{*} {\mathcal{\upomega}_{\mathbb{R}^4}}$ = $\mathcal{\upomega}_U$, with $\upomega_{\mathbb{R}^4}$ the standard symplectic form on $\mathbb{R}^4$(e.g. $dx_1 \wedge dy_1 + dx_2 \wedge dy_2)$ and $\pi: f^{*}C_{\mathbb{R}^4}^\infty \to C_U^\infty$.
\newline
\newline
This includes symplectic embeddings or immersions from U(2-dimensional manifold) into $\mathbb{R}^4$(4-dimensional manifold), paired with algebraic data.
\newline
\newline
Let's come back from the examples and talk about it again.
\newline
Verification: For finite limits in $\G$, the functor $\OX$ preserves them because mapping spaces and sheafification preserve limits.(Left exactness) For a cover $[(U_i,C_U^\infty) \to (M,C_M^\infty)]$ (where $[U_i]$ covers $M$), the sheaf condition holds: 
\begin{equation*}
    \OX(N,B)(M) \to \varprojlim \OX(N,B)(U_i)
\end{equation*}
is an equivalence, since $\OX(N,B)$ is a sheaf on $M$, and the topology on $\G$ aligns with the open cover topology on $M$.(Locality) For an admissible morphism $(i,\pi) : (U,C_U^\infty) \to (M,C_M^\infty)$(where $i: U \to M$ is an open embedding), $\OX(i,\pi)$ corresponds to restriction of sheaves which is an effective epimorphism in $Shv(M)$.
\newline

Note that: $\ind((\G)^{op})  \hookrightarrow \fun((\G),S)$. Since $\G$ has finite limits, representable functors are left exact, so S can be viewed in $\pro(\G)$.
\newline
For each $(N,B)$, $\OX(N,B)$ is a sheaf on M that describes symplectic maps from open subsets of M to N, along with compatible morphisms of $E_\infty$-rings from B to $C^{\infty}$.This construction is a precise realization of Lurie’s affine $\mathcal{G}$-scheme, tailored to the symplectic geometry context of $\G$.
\end{example}
\newline

\begin{defn}
An affine $\G$-scheme is a $\G$-structured topos ($\X,\mathcal{O}_{\X}$), if there exists an object $A \in \pro(\G$) such that:
\begin{equation*}
(\X,\mathcal{O}_{\X}) \cong \spec^{\G}(A).
\end{equation*}
The underlying $\infty$-topos $\X$ is the category of sheaves on $\G$ with respect to its Grothendieck topology, localized according to A. The structure sheaf $\mathcal{O}_{\X}$: $\G$ $\to$ $\X$ encodes the algebraic and symplectic data of A, with respect to the symplectic maps and $E_{\infty}$-ring morphisms defined in $\G$.
\end{defn}

Intuitively, an affine $\G$-scheme represents a space whose structure is globally determined by a single object in $\pro(\G)$, combining the symplectic geometry of a manifold with the algebraic structure of an $E_{\infty}$-ring sheaf.
\newline
\newline
We assumed that the reader is familiar with K-theory.
\begin{example}
K(A) be defined as the K-theory of the category of perfect A-modules, which are finitely presented modules in homotopical sense. (i.e. they behave like projective modules but in a derived setting). Let's take A from affine $\G$-scheme(from the definition above), where A is an $E_{\infty}$-ring equipped with symplectic data. The K-theory K(A) is computed from the perfect A-modules. This spectrum reflects algebraic properties of A, such as information about its module categories, enriched by the symplectic structure.

\end{example}

\begin{defn}
A $\G$-scheme is a $\G$-structured topos $(\X,\mathcal{O}_{\X})$ that is locally equivalent to affine $\G$-schemes. Specifically there exists a collection of open embeddings: 
\begin{equation*}
    (\X_i,\mathcal{O}_{\X_i}) \hookrightarrow (\X,\mathcal{O}_{\X}),
\end{equation*}
satisfying: Each $(\X_i,\mathcal{O}_{\X_i})$ is an affine $\G$-scheme, i.e. $(\X_i,\mathcal{O}_{\X_i}) \cong \spec^{\G}(A_i)$ for some $A_i \in \pro(\G)$, and the subtopoi $\X_i$ cover $\X$, meaning that the canonical map $\coprod \X_i \to \X$ is effective epimorphism.
\end{defn}
\newline
This definition implies that a $\G$-scheme can be glued together from the affine pieces, each carrying the symplectic and algebraic structure dictated by $\G$. The open embeddings reflect the topology generated by open symplectic embeddings in $\G$.
\newline
\begin{defn}
Consider a morphism f:($\X,\mathcal{O}_{\X}$) $\to$ ($\mathcal{Y},\mathcal{O}_{\mathcal{Y}}$) between $\G$-schemes, defined by a geometric morphism $f^*: \mathcal{Y} \to \X$ (with right adjoint $f_* $), a natural transformation $\beta : f^*\mathcal{O}_{\mathcal{Y}} \to \mathcal{O}_{\X}$, compatible with the symplectic structure.
\end{defn}
\newline

\begin{example}
For a $\G$-scheme $(\X,\mathcal{O}_{\X})$, consider an open subtopoi $U \subseteq \X$. The category of perfect modules over $\mathcal{O}_{\X}(U)$, denoted $\perf(\mathcal{O}_{\X}(U))$, is an $\infty$-category. The K-theory $K(\perf(\mathcal{O}_{\X}(U)))$ is a spectrum that captures the higher invariants of these modules, reflecting both their algebraic structure and the homotopical relationships between them. Since $\G$-schemes involve symplectic manifolds, we can refine $\perf(\mathcal{O}_{\X}(U))$ by considering subcategories where modules carry symplectic data(e.g. equipped with a symplectic forms or satisfying Lagrangian conditions). This defines a symplectic $\infty$-category, where objects are module with symplectic structures and morphisms preserve these structures up to higher homotopies. We will discuss it in the section[3].
\end{example}
\begin{defn}
A geometric morphism $f_*: \X \to \mathcal{Y}$ is closed embedding if $\X$ is equivalent to a closed subtopos of $\mathcal{Y}$.
\end{defn}
\begin{defn}
A morphism  (f, $\beta$):($\X, \mathcal{O}_X) \to (\mathcal{Y}, \mathcal{O}_{\mathcal{Y}}$) of $\G$-schemes is a closed immersion if it satisfies the following two conditions: 
\newline
\newline
(1): The geometric morphism $f_*: \X \to \mathcal{Y}$ is a closed embedding of $\infty$-topoi.
\newline
(2): For each object (M,A) $\in \G$, the component of the natural transformation $\beta_{(M,A)}$: $f^*\mathcal{O}_{\mathcal{Y}}(M,A) \to \OX(M,A)$ is an effective epimorphism in the $\infty$-topos $\X$. Here $f^*\mathcal{O}_{\mathcal{Y}}(M,A) \in \X$ is the pullback of the sheaf $\mathcal{O}_{\mathcal{Y}}(M,A) \in \mathcal{Y}$ along $f^*$ and $\OX(M,A) \in \X$ is the structure sheaf evaluated at (M,A). An effective epimorphism in an $\infty$-topos is a morphism that is surjective up to homotopy, meaning that $\OX(M,A)$ is a quotient of $f^*\mathcal{O}_{\mathcal{Y}}(M,A)$ in a categorical sense.
\newline
\newline
This condition ensures that $\OX$ is determined by $\mathcal{O}_{\mathcal{Y}}$ via a sheaf of ideals analogue, adapted to $E_{\infty}$-ring sheaves.
\end{defn}

\begin{defn}
A closed sub $\G$-scheme of a $\G$-scheme ($\mathcal{Y},\mathcal{O}_{\mathcal{Y}}$) is a $\G$-scheme ($\X, \mathcal{O}_X)$ together with a closed immersion (f, $\beta$): ($\X, \mathcal{O}_X) \to (\mathcal{Y},\mathcal{O}_{\mathcal{Y}}$).
\end{defn}
Since $\G$-schemes are glued from affine pieces, a closed subscheme can be constructed locally on affine charts and patched together, provided the closed embeddings and surjective sheaf maps are consistent across overlaps.

\section{Symplectic $\infty$-category}
\subsection {Perf($\X$)}
\begin{defn}
Define QCoh($\X$) as the $\infty$-category of quasi-coherent sheaves on ($\X,\OX$). Locally, on an affine open subtopos $Spec^{\G}(A_i) = (\X_i,\OX_i)$ where $A_i \in$ Pro($\G$), this corresponds to the category of modules over the local ring derived from $A_i$.
\end{defn}

\begin{defn}
Within QCoh($\X$), define Perf($\X$) as the full subcategory of perfect complexes. These are objects that are locally (in the $\infty$-topos sense) dualizable or equivalently, compact objects in QCoh($\X$). For an affine piece Spec$^{\G}$(A), Perf($\X$) restricts to Perf(A), the category of perfect A-modules.
\end{defn}

For an open subtopos U $\subseteq \X$, consider Perf($\OX(U))$. However, since $\OX: \G \to \X$ is a functor to the $\infty$-topos $\X$, $\OX(U)$ isn’t directly defined. Instead, interpret this via the restriction of $\OX$ to U, or consider U as corresponding to an open subset of a symplectic manifold M in the affine chart Spec$^{\G}$(A), where A = (M,$A_M$) and $A_M$ is a sheaf of $E_{\infty}$-rings(e.g. $C_M^{\infty}$). Locally, Perf($\OX(U))$) approximates Perf($A_{M}|_U$), the category of perfect modules over $A_{M}|_U$.
\newline
\newline
\begin{defn}
A $\G$-scheme ($\X,\OX$) is smooth over k if Perf($\X$) is smooth as a k-linear $\infty$-category.
\end{defn}
\begin{defn}
A $\G$-scheme ($\X,\OX$) is proper over k if Perf($\X$) is proper. (i.e. mapping spectra between compact objects are perfect over k.)
\end{defn}
\subsection{Symp($\X$) and K-theory of Symp($\X$)}
\begin{defn}
Category Symp$_d(\X$) with fixed d consists of the following data:
\newline
\newline
(1) Objects: Pairs (E,$\mathcal{\upomega}_E$), where $E \in Perf(\X)$ and $\mathcal{\upomega}_E: E \otimes E \to \OX[d]$ is non-degenerate skew-symmetric symplectic form, where $\OX[d]$ is the d-fold shift of the structure sheaf (viewed as an object in QCoh($\X$) via evaluation at a base object, or globally shifted) and d is the degree of the Calabi-Yau shift.
\newline
\newline
(2) Morphisms: Morphisms f: E $\to$ F in Perf($\X$) such that there exists a homotopy h: $\omega_F \circ $(f $\otimes$ f)  $\cong \mathcal{\omega}_E $, morphisms from a space with higher coherences, as in an $\infty$-category, meaning h is a path (in the homotopical sense) between these two maps from $E \otimes E$ to $\OX[d]$.
\newline
\newline
(3) Calabi-Yau Structure: Equip Symp$_d(\X$) with a d-Calabi-Yau Structure, where d arises from the symplectic form, ensuring $S_{symp} \cong [d]$.
\end{defn}
Let's define Symp($\X$). 
\begin{defn}
The graded category Symp($\X$) = $\bigoplus_{d}Symp_d(\X)$ is the direct sum (or coproduct) of stable $\infty$-categories Symp$_d(\X)$, indexed by d $\in \mathbb{Z}$:
\newline
Morphisms within Symp$_d(\X)$ are those that preserve the symplectic structure for that specific d. There are no morphisms between objects in Symp$_d(\X)$ and Symp$_{d'}(\X)$ when d $\neq d'$. The index d acts as a grading, separating the category into disjoint components based on the shift of the symplectic form.
\end{defn}
\begin{defn}
The K-theory of Symp$_d(\X$), denoted K(Symp$_d(\X$)), is the K-theory spectrum associated with the stable $\infty$-category Symp$_d(\X$).
\end{defn}
The homotopy groups of this spectrum give the K-theory groups: $\pi_n(K(Symp_d(\X)) = K_n(Symp_d(\X))$ for all n $\geq$ 0.
\newline
\newline

If $\OX = C_M^{\infty}$ and M is contractible, we get results similar(or equal) to previous Symplectic K-theory. We will discuss this(Symplectic K-theory) in the next paper. 
\newline
\newline
Thank you Anne Lee for providing an environment for a week of research. Thank you so much Taewan Kim for providing motivation and helping us with latex work. Note that: Taewan Kim will join the next paper !

Note that: This paper is just memos. We'll submit the completed paper next time. But we have other research projects(Check our page!), so we don't know when we're going to do this.

\bibliographystyle{unsrt}

\Large{Visiting our page! : \href {https://sites.google.com/view/pocariteikoku/home}{https://sites.google.com/view/pocariteikoku/home}}

\end{document}